\documentclass[12pt]{amsart}

\usepackage{enumitem}
\usepackage{psfrag}
\usepackage{graphicx}
\usepackage{pinlabel}
\usepackage{graphicx}
\usepackage{amsmath}
\usepackage{amssymb}
\usepackage{amscd}
\usepackage{autobreak}
\usepackage{picinpar}
\usepackage{times}
\usepackage{pb-diagram}
\usepackage{graphicx}
\usepackage{wrapfig}
\usepackage{xspace}
\usepackage{hyperref}
\usepackage{url}
\usepackage{caption}
\usepackage{subcaption}
\usepackage{fancyhdr}
\usepackage{overpic}
\usepackage{color}
\usepackage[square,comma,sort&compress,numbers]{natbib}
\usepackage{latexsym}
\usepackage{mathrsfs}
\usepackage{amsfonts}
\usepackage{hyperref}
\usepackage{amsthm}
\usepackage{appendix}
\usepackage{pdfsync}

\theoremstyle{definition}
\newtheorem{theorem}{Theorem}[section]
\newtheorem{lemma}[theorem]{Lemma}
\newtheorem{proposition}[theorem]{Proposition}
\newtheorem{corollary}[theorem]{Corollary}
\newtheorem{definition}[theorem]{Definition}

\newtheorem{remark}[theorem]{Remark}
\newtheorem*{theorem*}{Theorem}
\def\qed{\hfill{Q.E.D.}\smallskip}

\begin{document}

\title{\bf A discrete uniformization theorem for decorated piecewise Euclidean metrics on surfaces, II}
\author{Xu Xu, Chao Zheng}

\date{\today}

\address{School of Mathematics and Statistics, Wuhan University, Wuhan, 430072, P.R.China}
 \email{xuxu2@whu.edu.cn}

\address{School of Mathematics and Statistics, Wuhan University, Wuhan 430072, P.R. China}
\email{czheng@whu.edu.cn}



\begin{abstract}
In this paper, we study a natural discretization of the smooth Gaussian curvature on surfaces,
which is defined as the quotient of the angle defect and the area of a geodesic disk at a vertex of a polyhedral surface.
It is proved that each decorated piecewise Euclidean metric on surfaces with nonpositive Euler number is discrete conformal to a decorated piecewise Euclidean metric with this discrete curvature constant.
We further investigate the prescribing combinatorial curvature problem for a parametrization of this discrete curvature and prove some
Kazdan-Warner type results.
The main tools are Bobenko-Lutz's discrete conformal theory for decorated piecewise Euclidean metrics on surfaces and variational principles with constraints.
\end{abstract}

\maketitle

\section{Introduction}

The classical Gaussian curvature at a point $p$ in a Riemann surface can be defined as 
\begin{equation*}
R(p)=\lim_{r\rightarrow 0}\frac{12}{\pi r^4}(\pi r^2-A(r)),
\end{equation*}
where $A(r)$ is the area of the geodesic disk of radius $r$ at $p$.
If we apply this definition to a vertex $i$ of a piecewise Euclidean surface,
this gives a natural discretization of the classical Gaussian curvature (up to a constant)
\begin{equation}\label{Eq: curvature R}
R_i=\frac{K_i}{r_i^2},
\end{equation}
where $K_i=2\pi-\theta_i$ is the angle defect at $i$, $\theta_i$ is the cone angle at $i$ and $r_i$ is the radius of the geodesic disk at $i$.
We call $R_i$ as the discrete Gaussian curvature or combinatorial curvature.
It is natural to study the discrete uniformization theorem for the discrete Gaussian curvature $R_i$.
A good approach to this problem is working in the framework of decorated piecewise Euclidean metrics
recently introduced by Bobenko-Lutz \cite{BL}.

Suppose $S$ is a connected closed surface and $V$ is a finite non-empty subset of $S$,  the pair $(S, V)$ is called a marked surface.
A piecewise Euclidean metric (PE metric) $dist_{S}$ on the marked surface $(S,V)$
is a flat cone metric with the conic singularities contained in $V$.
A marked surface with a PE metric is called a PE surface, denoted by $(S,V, dist_{S})$.
The points in $V$ are called vertices of the PE surface.
A decoration $r$ on a PE surface $(S,V, dist_{S})$ is a choice of circle of radius $r_i$ at each vertex $i\in V$.
These circles in the decoration are called vertex-circles.
We denote a decorated PE surface by $(S,V, dist_{S}, r)$ and call the pair $(dist_S,r)$ a decorated PE metric on the marked surface $(S,V)$.
In this paper, we focus on the case that each pair of vertex-circles is separated.

\begin{theorem}\label{Thm: existence 1}
Let $(dist_S,r)$ be a decorated PE metric on a marked surface $(S,V)$ with Euler number $\chi(S)\leq0$.
Let $\overline{R}\leq 0$ be a function defined on $V$ satisfying $\overline{R}\not\equiv0$ if $\chi(S)<0$ and
$\overline{R}\equiv0$ if $\chi(S)=0$.
Then there exists a unique discrete conformal equivalent decorated PE metric $(\widetilde{dist_S},\widetilde{r})$
on $(S,V)$ with the prescribed discrete Gaussian curvature $\overline{R}$ (up to scaling if $\chi(S)=0$).
\end{theorem}

Theorem \ref{Thm: existence 1} is a discrete analogue of Kazdan-Warner's results in \cite{KW1,KW2}.
As a corollary of Theorem \ref{Thm: existence 1},
we have the following discrete uniformization theorem for the discrete Gaussian curvature $R_i$ on decorated PE surfaces.
\begin{corollary}\label{Cor: existence 1}
For any decorated PE metric $(dist_S,r)$ on a marked surface $(S,V)$ with Euler number $\chi(S)\leq0$,
there exists a unique discrete conformal equivalent decorated PE metric $(\widetilde{dist_S},\widetilde{r})$ on $(S,V)$
with constant discrete Gaussian curvature $R$ (up to scaling if $\chi(S)=0$).
\end{corollary}

The combinatorial curvature $R$ in (\ref{Eq: curvature R}) was first introduced by Ge-Xu \cite{GX2} for Thurston's circle packing metrics on triangulated surfaces.
After that, there are lots of research activities on the combinatorial curvature $R$ on surfaces. See \cite{GJ, GX2,GX IMRN,GX JFA, Xu1,XZ TAMS, XZ CVPDE} and others for example.
Most of these works gave equivalent conditions for the existence of discrete conformal metrics with prescribed
combinatorial curvature $R$ via combinatorial curvature flows.
In Theorem \ref{Thm: existence 1} and Corollary \ref{Cor: existence 1}, we give some sufficient conditions
for the existence involving only the prescribed combinatorial curvatures and the topology of the surfaces.

Following \cite{GX2}, we further introduce the following
parameterized combinatorial curvature for the decorated PE metrics on surfaces
\begin{equation}\label{Eq: curvature R_alpha}
R_{\alpha, i}=\frac{K_i}{r_i^\alpha},
\end{equation}
where $\alpha\in \mathbb{R}$ is a constant.
If $\alpha=2$, then $R_{2, i}$ is the combinatorial curvature $R_i$ defined by (\ref{Eq: curvature R}).
We call $R_{\alpha}$ as the combinatorial $\alpha$-curvature.


\begin{theorem}\label{Thm: existence 2}
Let $(dist_S,r)$ be a decorated PE metric on a marked surface $(S,V)$.
Suppose $\alpha\in \mathbb{R}$ is a constant and $\overline{R}: V\rightarrow\mathbb{R}$ is a given function.
There exists a discrete conformal equivalent decorated PE metric $(\widetilde{dist_S},\widetilde{r})$ with combinatorial $\alpha$-curvature $\overline{R}$ if one of the following conditions is satisfied
\begin{description}
  \item[(1)] $\chi(S)>0,\ \alpha<0,\  \overline{R}>0$;
  \item[(2)] $\chi(S)<0,\ \alpha\neq0,\  \overline{R}\leq 0,\ \overline{R}\not\equiv 0$;
  \item[(3)] $\chi(S)=0,\ \alpha\neq0,\  \overline{R}\equiv0$;
  \item[(4)] $\alpha=0$, $\overline{R}\in (-\infty, 2\pi)$, $\sum_{i\in V}\overline{R}_{i}=2\pi \chi(S)$.
\end{description}
If $\alpha\overline{R}\leq0$,  the decorated PE metric $(\widetilde{dist_S},\widetilde{r})$ is unique (up to scaling if $\alpha\overline{R}\equiv 0$).
\end{theorem}

Theorem \ref{Thm: existence 2} is a generalization of Theorem \ref{Thm: existence 1}. Specially, if $\alpha=2$, then the cases \textbf{(2)}
and \textbf{(3)} in Theorem \ref{Thm: existence 2} are reduced to Theorem \ref{Thm: existence 1}.
By the relationship of the combinatorial $\alpha$-curvature $R_\alpha$ and the angle defect $K$,
the cases \textbf{(3)} and \textbf{(4)} in Theorem \ref{Thm: existence 2} are covered by Bobenko-Lutz's work \cite{BL}.
In the following, we just prove the cases \textbf{(1)} and \textbf{(2)} of Theorem \ref{Thm: existence 2}.

\begin{remark}
Since Bobenko-Lutz's discrete conformal theory of decorated PE metrics also applies to Luo's vertex scalings and thus generalizes Gu-Luo-Sun-Wu's results in \cite{Gu1} and Springborn's results in \cite{Springborn},
Theorem \ref{Thm: existence 2} also generalizes the authors' results in \cite{XZ CVPDE}.
\end{remark}

As a corollary of Theorem \ref{Thm: existence 2},
we have the following discrete uniformization theorem for the combinatorial $\alpha$-curvature $R_\alpha$.

\begin{corollary}\label{Cor: existence 2}
Suppose $(S,V)$ is a marked surface with a decorated PE metric $(dist_S,r)$ and $\alpha\in \mathbb{R}$ is a constant.
\begin{description}
\item[(1)]
If $\alpha\chi(S)\leq0$, there exists a unique discrete conformal equivalent decorated PE metric $(\widetilde{dist_S},\widetilde{r})$ with constant combinatorial $\alpha$-curvature $R_\alpha$ (up to scaling if $\alpha\chi(S)=0$).
\item[(2)]
If $\alpha<0$ and $\chi(S)<0$, there exists a discrete conformal equivalent decorated PE metric $(\widetilde{dist_S},\widetilde{r})$ with negative constant combinatorial $\alpha$-curvature $R_\alpha$.
\end{description}
\end{corollary}

The paper is organized as follows.
In Section \ref{Sec: DC theory}, we briefly recall Bobenko-Lutz's discrete conformal theory for the decorated PE metrics on surfaces.
Then we prove the global rigidity of decorated PE metrics with respect to the combinatorial $\alpha$-curvature on a marked surface.
In Section \ref{Sec: existence},
we first deform the combinatorial $\alpha$-curvature $R_\alpha$ in (\ref{Eq: curvature R_alpha}) and give Theorem \ref{Thm: existence 3}, which is equivalent to Theorem \ref{Thm: existence 2}.
Then we translate Theorem \ref{Thm: existence 3} into an optimization problem with constraints, i.e., Theorem \ref{Thm: inequality constraints}.
Using a classical result from calculus, i.e., Theorem \ref{Thm: calculus}, we translate Theorem \ref{Thm: inequality constraints} into Theorem \ref{Thm: main 2}.
In the end, with the help of the asymptotical expression for the energy function $\mathcal{E}$ in Lemma \ref{Lem: E decomposition}  obtained by the authors in \cite{XZ}, we prove Theorem \ref{Thm: main 2}.
\\
\\
\textbf{Acknowledgements}\\[8pt]
The first author thanks Professor Feng Luo for his invitation to the workshop
``Discrete and Computational Geometry, Shape Analysis, and Applications" taking place
at Rutgers University, New Brunswick from May 19th to May 21st, 2023.
The first author also thanks Carl O. R. Lutz for helpful communications during the workshop.

\section{Rigidity of decorated PE metrics}
\label{Sec: DC theory}

\subsection{Discrete conformal equivalence and Bobenko-Lutz's discrete conformal theory}


Let $\mathcal{T}={(V,E,F)}$ be a triangulation of a marked surface $(S, V)$,
where $V,E,F$ represent the sets of vertices, edges and faces respectively.
A triangulation $\mathcal{T}$ for a PE surface $(S,V, dist_S)$ is a geodesic triangulation if the edges are geodesics in the PE metric $dist_S$.
We use one index to denote a vertex (such as $i$), two indices to denote an edge (such as $\{ij\}$) and three indices to denote a face (such as $\{ijk\}$) in the triangulation $\mathcal{T}$.
The PE metric $dist_{S}$ on a PE surface with a geodesic triangulation $\mathcal{T}$ defines a length map $l: E\rightarrow \mathbb{R}_{>0}$ such that $l_{ij}, l_{ik}, l_{jk}$ satisfy the triangle inequalities for any triangle $\{ijk\}\in F$.
Conversely, given a function $l: E\rightarrow \mathbb{R}_{>0}$ satisfying the triangle inequalities for any face $\{ijk\}\in F$,
one can construct a PE metric on a triangulated surface by isometrically gluing Euclidean triangles along edges in pairs.
In the following, we use $l: E\rightarrow \mathbb{R}_{>0}$ to denote a PE metric and use $(l,r)$ to denote a decorated PE metric on a triangulated surface $(S,V,\mathcal{T})$.

\begin{definition}[\cite{BL}, Proposition 2.2]
\label{Def: DCE}
Let $\mathcal{T}$ be a triangulation of a marked surface $(S,V)$.
Two decorated PE metrics $(l,r)$ and $(\widetilde{l},\widetilde{r})$ on $(S,V, \mathcal{T})$ are discrete conformal equivalent
if and only if there exists a discrete conformal factor $u\in \mathbb{R}^V$ such that
\begin{equation}\label{Eq: DCE1}
\widetilde{r}_i=e^{u_i}r_i,
\end{equation}
\begin{equation}\label{Eq: DCE2}
\widetilde{l}_{ij}^2
=(e^{2u_i}-e^{u_i+u_j})r^2_i
+(e^{2u_j}-e^{u_i+u_j})r^2_j
+e^{u_i+u_j}l_{ij}^2
\end{equation}
for any edge $\{ij\}\in E$.
\end{definition}

\begin{remark}\label{Rmk: inversive distance}
For any two circles $C_i$ and $C_j$ in the Euclidean plane, one can define the inversive distance
$I_{ij}=\frac{l^2_{ij}-r^2_i-r^2_j}{2r_ir_j}$,
where $l_{ij}$ is the distance of the centers of the two circles and $r_i$, $r_j$ are the radii of $C_i, C_j$ respectively.
The inversive distance is invariant under M\"{o}bius transformations \cite{Coxeter}.
Denote the inversive distance of two vertex-circles in $(l,r)$ and $(\widetilde{l},\widetilde{r})$ as $I$ and $\widetilde{I}$ respectively.
If $(l,r)$ and $(\widetilde{l},\widetilde{r})$ are discrete conformal equivalent in the sense of Definition \ref{Def: DCE},
it is shown \cite{BL} that $I=\widetilde{I}$.
Since each pair of vertex-circles is required to be separated, it is easy to see that $I>1$.
Therefore, the discrete conformal equivalent decorated PE metrics on triangulated surfaces in Definition \ref{Def: DCE} can be taken as the separated inversive distance circle packing metrics introduced by Bowers-Stephenson \cite{BS}.
Please refer to \cite{CLXZ,Guo,Luo3,Xu AIM,Xu MRL} for more properties of the inversive distance circle packing metrics on triangulated surfaces.
\end{remark}

For any decorated triangle $\{ijk\}$,
there is a unique circle $C_{ijk}$ simultaneously orthogonal to the three vertex-circles at the vertices $i,j,k$ \cite{Glickenstein preprint}.
This circle $C_{ijk}$ is called as the face-circle of the decorated triangle $\{ijk\}$.
Denote $\alpha_{ij}^k$ as the interior intersection angle of the face-circle $C_{ijk}$ and the edge $\{ij\}$.
The edge $\{ij\}$, shared by two adjacent decorated triangles $\{ijk\}$ and $\{ijl\}$, is called weighted Delaunay if
\begin{equation*}
\alpha_{ij}^k+\alpha_{ij}^l\leq \pi.
\end{equation*}
The triangulation $\mathcal{T}$ is called weighted Delaunay in the decorated PE metric $(dist_S,r)$ if every edge in the triangulation is weighted Delaunay.
Here we take the definition of weighted Delaunay triangulation from Bobenko-Lutz \cite{BL}.
There are other equivalent definitions for the weighted Delaunay triangulation using the signed distance of the center of $C_{ijk}$ to the edges.
Please refer to  \cite{CLXZ, Glickenstein DCG,Glickenstein JDG,Glickenstein preprint,GT} and others.

Note that the combinatorial $\alpha$-curvature $R_\alpha$ in (\ref{Eq: curvature R_alpha}) is independent of the geodesic triangulations of a decorated PE surface $(S,V, dist_{S}, r)$.
In general, the existence of decorated PE metrics with prescribed combinatorial $\alpha$-curvatures on triangulated surfaces can not be guaranteed if the triangulation is fixed.
In the following, we work with a generalization of the discrete conformal equivalence in Definition \ref{Def: DCE},
introduced by Bobenko-Lutz \cite{BL},
which allows the triangulation of the marked surface to be changed under the weighted Delaunay condition.

\begin{definition}[\cite{BL}, Definition 4.11]
\label{Def: GDCE}
Two decorated PE metrics $(dist_{S},r)$ and $(\widetilde{dist}_{S},\widetilde{r})$ on the marked surface $(S,V)$ are discrete conformal equivalent if there is a sequence of triangulated decorated PE surfaces
$(\mathcal{T}^0,l^0,r^0),...,(\mathcal{T}^N,l^N,r^N)$ such that
\begin{description}
  \item[(1)] the decorated PE metric of $(\mathcal{T}^0,l^0,r^0)$ is $(dist_{S},r)$ and the decorated PE metric of $(\mathcal{T}^N,l^N,r^N)$ is $(\widetilde{dist}_{S},\widetilde{r})$,
  \item[(2)] each $\mathcal{T}^n$ is a weighted Delaunay triangulation of the decorated PE surface $(\mathcal{T}^n,l^n,r^n)$,
  \item[(3)] if $\mathcal{T}^n=\mathcal{T}^{n+1}$, then there is a discrete conformal factor $u\in \mathbb{R}^V$ such that $(\mathcal{T}^n,l^n,r^n)$ and $(\mathcal{T}^{n+1},l^{n+1},r^{n+1})$ are related by (\ref{Eq: DCE1}) and (\ref{Eq: DCE2}),
  \item[(4)] if $\mathcal{T}^n\neq\mathcal{T}^{n+1}$, then $\mathcal{T}^n$ and $\mathcal{T}^{n+1}$ are two different weighted Delaunay triangulations of the same decorated PE surface.
\end{description}
\end{definition}

Definition \ref{Def: GDCE} gives an equivalence relationship for decorated PE metrics on a marked surface.
The equivalence class of a decorated PE metric $(dist_S,r)$ on $(S,V)$ is also called as the discrete conformal class of $(dist_S,r)$ and denoted by $\mathcal{D}(dist_S,r)$.

\begin{lemma}[\cite{BL}]
The discrete conformal class $\mathcal{D}(dist_S,r)$ of a decorated PE metric $(dist_S,r)$ on the marked surface $(S,V)$ is parameterized by $\mathbb{R}^V=\{u: V\rightarrow \mathbb{R}\}$.
\end{lemma}
For simplicity, for any $(\widetilde{dist}_S,\widetilde{r})\in \mathcal{D}(dist_S,r)$,
we denote it by $(dist_S(u),r(u))$ for some $u\in \mathbb{R}^V$.
Set
\begin{equation*}
\mathcal{C}_\mathcal{T}(dist_{S},r)
=\{u\in \mathbb{R}^V |\ \mathcal{T}\ \text{is a weighted Delaunay triangulation of}\ (S,V,dist_S(u),r(u))\}.
\end{equation*}

\begin{lemma}[\cite{BL}]\label{Lem: finite decomposition}
The set
$$J=\{\mathcal{T}| \mathcal{C}_{\mathcal{T}}(dist_{S},r)\ \text{has non-empty interior in}\ \mathbb{R}^V\}$$
is a finite set, $\mathbb{R}^V=\cup_{\mathcal{T}_i\in J}\mathcal{C}_{\mathcal{T}_i}(dist_{S},r)$ and
each $\mathcal{C}_{\mathcal{T}_i}(dist_{S},r)$ is homeomorphic to a polyhedral cone (with its apex removed)
and its interior is homeomorphic to $\mathbb{R}^V$.
\end{lemma}

\subsection{The extended energy function}

There exist geometric relationships between the decorated triangles and $3$-dimensional generalized hyperbolic polyhedra.
Specially, there is a generalized hyperbolic tetrahedra in $\mathbb{H}^3$ with one ideal vertex and three hyper-ideal vertices corresponding to a decorated triangle $\{ijk\}$.
Denote $\mathrm{Vol}(ijk)$ as the truncated volume of this generalized hyperbolic tetrahedra.
The truncated volume $\mathrm{Vol}(ijk)$ can be characterized by an explicit formula.
Please refer to \cite{BL, Sp1} for more details.

Set
\begin{equation*}
\begin{aligned}
F_{ijk}(u_i,u_j,u_k)
=&-2\mathrm{Vol}(ijk)+\theta_{jk}^iu_i+\theta_{ki}^ju_j
+\theta_{ij}^ku_k\\
&+(\frac{\pi}{2}-\alpha_{ij}^k)\lambda_{ij}
+(\frac{\pi}{2}-\alpha_{ki}^j)\lambda_{ki}
+(\frac{\pi}{2}-\alpha_{jk}^i)\lambda_{jk},
\end{aligned}
\end{equation*}
where $\theta_{jk}^i$ is the inner angle of the decorated triangle $\{ijk\}$ at the vertex $i$ and $\lambda_{ij}=\cosh^{-1} I_{ij}$.
By the Schl\"{a}fli formula, we have
$$\nabla F_{ijk}=(\theta_{jk}^i,\theta_{ki}^j, \theta_{ij}^k)$$
and
$$F_{ijk}((u_i,u_j,u_k)+c(1,1,1))=F_{ijk}(u_i,u_j,u_k)+c\pi$$
for $c\in \mathbb{R}$.
On a decorated PE surface $(S,V,l,r)$ with a weighted Delaunay triangulation $\mathcal{T}$,
Bobenko-Lutz \cite{BL} defined the following function
\begin{equation}\label{Eq: F1}
\mathcal{H}_{\mathcal{T}}(u)
=\sum_{\{ijk\}\in F}F_{ijk}(u_i,u_j,u_k)
=-2\sum_{\{ijk\}\in F}\mathrm{Vol}(ijk)+\sum_{i\in V}\theta_iu_i+\sum_{\{ij\}\in E}(\pi-\alpha_{ij})\lambda_{ij},
\end{equation}
where $\theta_i=\sum_{\{ijk\}\in F}\theta^i_{jk}$ and $\alpha_{ij}=\alpha_{ij}^k+\alpha_{ij}^l$.
It should be mentioned that the function $\mathcal{H}_{\mathcal{T}}(u)$ in (\ref{Eq: F1}) differs from its original definition in \cite{BL} (Equation 4-9) by some constant.
Then
\begin{equation}\label{Eq: property of H ijk}
\mathcal{H}_{\mathcal{T}}(u+c\mathbf{1})
=\mathcal{H}_{\mathcal{T}}(u)+c|F|\pi
\end{equation}
for $c\in \mathbb{R}$.
By the definition of $\mathcal{H}_{\mathcal{T}}$,
the following energy function
\begin{equation*}
\mathcal{E}_{\mathcal{T}}(u)
=-\mathcal{H}_{\mathcal{T}}(u)+2\pi\sum_{i\in V}u_i
\end{equation*}
is well-defined on $\mathcal{C}_\mathcal{T}(dist_{S},r)$ with
$\nabla_{u_i} \mathcal{E}_{\mathcal{T}}
=2\pi-\theta_i=K_i$.
Moreover,
\begin{equation}\label{Eq: property of E}
\mathcal{E}_{\mathcal{T}}(u+c\mathbf{1})
=\mathcal{E}_{\mathcal{T}}(u)+2c\pi\chi(S)
\end{equation}
for $c\in \mathbb{R}$.

\begin{theorem}[\cite{BL}, Proposition 4.13]
\label{Thm: extended H}
For a discrete conformal factor $u\in \mathbb{R}^V$, let $\mathcal{T}$
be a weighted Delaunay triangulation of the decorated PE surface $(S,V,dist_S(u),r(u))$.
The map
\begin{equation}\label{Eq: extended H}
\mathcal{H} :\  \mathbb{R}^V\rightarrow \mathbb{R},\ \quad
u\mapsto \mathcal{H}_{\mathcal{T}}(u)
\end{equation}
is well-defined, concave, and twice continuously differentiable over $\mathbb{R}^V$.
\end{theorem}
Therefore, the function $\mathcal{E}_{\mathcal{T}}(u)$ defined on $\mathcal{C}_\mathcal{T}(dist_{S},r)$ can be extended to be $\mathcal{E}(u)$ defined on $\mathbb{R}^V$ by the following formula
\begin{equation}\label{Eq: extended E}
\mathcal{E}(u)
=-\mathcal{H}(u)+2\pi\sum_{i\in V}u_i.
\end{equation}

\subsection{Rigidity of decorated PE metrics}

A basic problem on the combinatorial $\alpha$-curvature is to understand the relationships between the decorated PE metrics and the combinatorial $\alpha$-curvatures.
The following theorem shows the global rigidity of decorated PE metrics with respect to the combinatorial $\alpha$-curvature on a marked surface,
which corresponds to the rigidity parts of Theorem \ref{Thm: existence 2}.

\begin{theorem}\label{Thm: rigidity}
Suppose $(S,V)$ is a marked surface with a decorated PE metric $(dist_S,r)$, $\alpha\in \mathbb{R}$ is a constant and $\overline{R}: V\rightarrow\mathbb{R}$ is a given function.
\begin{description}
\item[(1)]
If $\alpha\overline{R}\equiv 0$, then there exists at most one discrete conformal factor $u^*\in \mathbb{R}^V$ up to scaling such that the decorated PE metric $(dist_S(u^*),r(u^*))$ in the discrete conformal class $\mathcal{D}(dist_S,r)$
 has the combinatorial $\alpha$-curvature $\overline{R}$.
\item[(2)]
If $\alpha\overline{R}\leq 0$ and $\alpha\overline{R}\not\equiv 0$, then there exists at most one discrete conformal factor $u^*\in \mathbb{R}^V$
such that the decorated PE metric $(dist_S(u^*),r(u^*))$ in the discrete conformal class $\mathcal{D}(dist_S,r)$ has the combinatorial $\alpha$-curvature $\overline{R}$.
\end{description}
\end{theorem}
\proof
By Theorem \ref{Thm: extended H}, the following function
\begin{equation}
\mathbb{E}(u)
=-\mathcal{H}(u)+\int^u_{u_0}
\sum_{i\in V}(2\pi-\overline{R}_ir_i^\alpha)du_i
\end{equation}
is well-defined and twice continuously differentiable over $\mathbb{R}^V$,
where $r_i=e^{u_i}r^0_i$ and $r^0$ is the initial data.
By direct calculations, we have
\begin{equation*}
\nabla_{u_i}\mathbb{E}
=-\sum_{\{ijk\}\in F}\theta^i_{jk}
+2\pi-\overline{R}_ir_i^\alpha
=K_i-\overline{R}_ir_i^\alpha.
\end{equation*}
Therefore, for $u^*\in \mathbb{R}^V$, the decorated PE metric in $(dist_S(u^*),r(u^*))$ has the combinatorial $\alpha$-curvature $\overline{R}$ if and only if $\nabla_{u_i}\mathbb{E}(u^*)=0, \forall i\in V$.
Moreover,
\begin{equation*}
\mathrm{Hess}_u\ \mathbb{E}
=-\mathrm{Hess}_u\ \mathcal{H} -\alpha \left(
 \begin{array}{ccc}
  \overline{R}_1r^\alpha_1  &  &\\
     & \ddots &  \\
      & & \overline{R}_{|V|}r^\alpha_{|V|} \\
  \end{array}
 \right).
\end{equation*}
The equality (\ref{Eq: property of H ijk}) and Theorem \ref{Thm: extended H} imply that $\mathrm{Hess}_u\mathcal{H}\leq0$ with kernel $\{c\mathbf{1}^\mathrm{T}\in \mathbb{R}^V|c\in \mathbb{R}\}$.
If $\alpha\overline{R}\equiv 0$, then $\mathrm{Hess}_u\ \mathbb{E}$ is positive semi-definite with kernel $\{c\mathbf{1}^\mathrm{T}\in \mathbb{R}^V|c\in \mathbb{R}\}$ and hence $\mathbb{E}$ is convex on $\mathbb{R}^V$ and strictly convex on $\{\sum_{i\in V}u_i=0\}$.
If $\alpha\overline{R}\leq 0$ and $\alpha\overline{R}\not\equiv 0$,
then $\mathrm{Hess}_u\ \mathbb{E}$ is positive definite and hence $\mathbb{E}$ is strictly convex on $\mathbb{R}^V$.
The conclusion follows from the following result from calculus.

\noindent\textbf{Lemma:}
If $f:\Omega \rightarrow \mathbb{R}$ is a $C^1$-smooth  strictly convex function on an open convex set $\Omega \subset \mathbb{R}^n$,
then its gradient $\nabla f:\Omega \rightarrow \mathbb{R}^n$ is injective.
\qed

\begin{remark}
For a decorated PE surfaces $(S,V,l,r)$ with a fixed triangulation $\mathcal{T}$,
the global rigidity of the inversive distance circle packing metrics with respect to the combinatorial $\alpha$-curvature $R_\alpha$ has been proved by Ge-Jiang \cite{GJ} and Ge-Xu \cite{GX2}.
They extended the function $\mathbb{E}(u)$ by extending the inner angles of a triangle by constants.
This approach was introduced by Bobenko-Pinkall-Springborn \cite{BPS} for Luo's vertex scalings and further
developed by Luo \cite{Luo3} for Bowers-Stephenson's inversive distance circle packings.
Here we take another approach introduced by Bobenko-Lutz \cite{BL}
to extend the function $\mathbb{E}(u)$, in which we change the triangulation of the marked surface under the weighted Delaunay condition.
This approach was first introduced by Gu-Luo-Sun-Wu \cite{Gu1} and Gu-Guo-Luo-Sun-Wu \cite{Gu2}
to establish the discrete uniformization theorem for Luo's vertex scalings on surfaces.
The first approach can not ensure the triangles being non-degenerate, while the second approach can.
\end{remark}

\section{Existence of decorated PE metrics}
\label{Sec: existence}

\subsection{Variational principles with constraints}

In this subsection, to simplify the calculations, we deform the combinatorial $\alpha$-curvature $R_\alpha$ in (\ref{Eq: curvature R_alpha}) and give Theorem \ref{Thm: existence 3} which is equivalent to Theorem \ref{Thm: existence 2}.
Then we translate Theorem \ref{Thm: existence 3} into an optimization problem with inequality constraints by variational principles, which involves the function $\mathcal{E}(u)$ defined in (\ref{Eq: extended E}).

For an initial decorated PE metric $(l^0,r^0)$,
the combinatorial $\alpha$-curvature is $K^0_i/(r^0_i)^\alpha$.
Suppose a decorated PE metric $(l,r)$ is discrete conformal equivalent to $(l^0,r^0)$,
then $r_i=e^{u_i}r_i^0$ by (\ref{Eq: DCE1}).
The combinatorial $\alpha$-curvature of the decorated PE metric $(l,r)$ can be written as
\begin{equation*}
R_{\alpha,i}
=\frac{K_i}{r^\alpha_i}
=\frac{K_i}{(r^0_i)^\alpha e^{\alpha u_i}}.
\end{equation*}
For simplicity, set
$$\mathcal{R}_{\alpha,i}=R_{\alpha,i}(r^0_i)^\alpha.$$
Then
\begin{equation}\label{Eq: key 4}
\mathcal{R}_{\alpha,i}=\frac{K_i}{e^{\alpha u_i}}.
\end{equation}
We also call $\mathcal{R}_{\alpha}$ as the combinatorial $\alpha$-curvature.
Note that $(r^0_i)^\alpha>0$, then the signs of $\mathcal{R}_{\alpha,i}$ and $R_{\alpha,i}$ are the same for any $i\in V$.
Denote $\overline{\mathcal{R}}$ as the prescribed combinatorial $\alpha$-curvature corresponding to $\mathcal{R}_{\alpha}$.
Then $\overline{\mathcal{R}}_i=\overline{R}_i(r^0_i)^\alpha$ and the signs of $\overline{\mathcal{R}}_i$ and $\overline{R}_i$ are the same.
Hence, to prove Theorem \ref{Thm: existence 2}, we just need to prove the following theorem.
\begin{theorem}\label{Thm: existence 3}
For any decorated PE metric $(dist_S,r)$ on a marked surface $(S,V)$,
there is a discrete conformal equivalent decorated PE metric $(\widetilde{dist_S},\widetilde{r})$ with combinatorial $\alpha$-curvature $\overline{\mathcal{R}}$ if one of the following conditions is satisfied
\begin{description}
\item[(1)] $\chi(S)>0,\ \alpha<0,\  \overline{\mathcal{R}}>0$;
\item[(2)] $\chi(S)<0,\ \alpha\neq0,\  \overline{\mathcal{R}}\leq 0,\ \overline{\mathcal{R}}\not\equiv 0$.
\end{description}
\end{theorem}

Since the angle defect $K$ satisfies the following discrete Gauss-Bonnet formula (\cite{Chow-Luo}, Proposition 3.1)
\begin{equation}\label{Eq: Gauss-Bonnet}
\sum_{i\in V}K_i=2\pi\chi(S),
\end{equation}
then the combinatorial $\alpha$-curvature $\mathcal{R}_\alpha$ in (\ref{Eq: key 4}) satisfies the following discrete Gauss-Bonnet formula
\begin{equation*}
\sum_{i\in V}\mathcal{R}_ie^{\alpha u_i}=2\pi\chi(S).
\end{equation*}
Therefore, if $\overline{\mathcal{R}}\in \mathbb{R}^V$ is the combinatorial $\alpha$-curvature of some decorated PE metric discrete conformal to $(l,r)$ on $(S,V)$, then
\begin{equation*}
\sum_{i\in V}\overline{\mathcal{R}}_ie^{\alpha u_i}=2\pi\chi(S).
\end{equation*}
Let $\alpha\in \mathbb{R}$ be a non-zero constant.
Set
\begin{equation}\label{A}
\mathcal{A}=\{u\in \mathbb{R}^V|0>\sum_{i\in V}\overline{\mathcal{R}}_i e^{\alpha u_i}\geq 2\pi\chi(S),\  \overline{\mathcal{R}}\leq0,\ \overline{\mathcal{R}}\not\equiv0\},
\end{equation}
\begin{equation}\label{B}
\mathcal{B}=\{u\in \mathbb{R}^V|0<\sum_{i\in V}\overline{\mathcal{R}}_i e^{\alpha u_i}\leq 2\pi\chi(S),\ \overline{\mathcal{R}}>0\},
\end{equation}
\begin{equation}\label{C}
\mathcal{C}=\{u\in \mathbb{R}^V|\sum_{i\in V} \overline{\mathcal{R}}_i e^{\alpha u_i}\leq 2\pi\chi(S)<0,\ \overline{\mathcal{R}}\leq0,\ \overline{\mathcal{R}}\not\equiv0 \}.
\end{equation}

\begin{proposition}\label{Prop: ABC}
The sets $\mathcal{A},\ \mathcal{B}$ and $\mathcal{C}$ are unbounded closed subsets of $\mathbb{R}^V$.
\end{proposition}
\proof
We only prove this proposition for the set $\mathcal{A}$ and the proofs for the sets $\mathcal{B}$ and $\mathcal{C}$ are similar.

\noindent\textbf{(I):}
To prove the closeness of the set $\mathcal{A}$ in $\mathbb{R}^V$, we just need to show $\mathcal{A}=\overline{\mathcal{A}}$,
where $\overline{\mathcal{A}}$ represents the closure of the set $\mathcal{A}$ in $\mathbb{R}^V$.
Suppose $\{u_{i,n}\}_{n\in \mathbb{N}}$ is a sequence in $\mathcal{A}$ such that $\lim_{n\rightarrow +\infty}u_{i,n}=\lambda_i\in \mathbb{R}, \forall i\in V$.
It is direct to see that
$\lim_{n\rightarrow +\infty} \sum_{i\in V}\overline{\mathcal{R}}_i e^{\alpha u_{i,n}}
=\sum_{i\in V}\overline{\mathcal{R}}_i e^{\alpha \lambda_i}\geq 2\pi\chi(S).$
Note that the definition of $\mathcal{A}$ in (\ref{A}) shows $\overline{\mathcal{R}}\leq0,\ \overline{\mathcal{R}}\not\equiv0$.
This implies that there exists $i_0\in V$ such that $\overline{\mathcal{R}}_{i_0}<0$.
Then
\begin{equation*}
\lim_{n\rightarrow +\infty} \sum_{i\in V}\overline{\mathcal{R}}_i e^{\alpha u_{i,n}}
=\sum_{i\in V}\overline{\mathcal{R}}_i e^{\alpha \lambda_{i}}\leq\overline{\mathcal{R}}_{i_0} e^{\alpha \lambda_{i_0}} <0.
\end{equation*}
This implies $\lambda=(\lambda_1,...,\lambda_{|V|})\in \mathcal{A}$ and hence $\mathcal{A}=\overline{\mathcal{A}}$.
Therefore, the set $\mathcal{A}$ is a closed subset of $\mathbb{R}^V$.

\noindent\textbf{(II):}
If $u\in \mathcal{A}$, for any $c\in \mathbb{R}$, we have
\begin{equation*}
\sum_{i\in V} \overline{\mathcal{R}}_i e^{\alpha (u_i+c)}=e^{\alpha c}\sum_{i\in V} \overline{\mathcal{R}}_i e^{\alpha u_i}<0.
\end{equation*}
If $\alpha<0$, $u\in \mathcal{A}$, then
\begin{equation*}
\sum_{i\in V} \overline{\mathcal{R}}_i e^{\alpha (u_i+c)}=e^{\alpha c}\sum_{i\in V} \overline{\mathcal{R}}_i e^{\alpha u_i}\geq 2\pi\chi(S)
\end{equation*}
is equivalent to
\begin{equation*}
c\geq\frac{1}{\alpha}\log\frac{2\pi\chi(S)} {\sum_{i\in V} \overline{\mathcal{R}}_i e^{\alpha u_i}}.
\end{equation*}
This implies that the ray $\{u+c\mathbf{1}|c\geq\frac{1}{\alpha}\log\frac{2\pi\chi(S)} {\sum_{i\in V} \overline{\mathcal{R}}_i e^{\alpha u_i}},\ \alpha<0 \}$ stays in the set $\mathcal{A}$.
Hence, the set $\mathcal{A}$ is unbounded  if $\alpha<0$.

If $\alpha>0$, for $u\in \mathcal{A}$, we have
\begin{equation*}
\sum_{i\in V} \overline{\mathcal{R}}_i e^{\alpha (u_i+c)}=e^{\alpha c}\sum_{i\in V} \overline{\mathcal{R}}_i e^{\alpha u_i}\geq 2\pi\chi(S)
\end{equation*}
is equivalent to
\begin{equation*}
c\leq\frac{1}{\alpha}\log\frac{2\pi\chi(S)} {\sum_{i\in V} \overline{\mathcal{R}}_i e^{\alpha u_i}}.
\end{equation*}
This implies that the ray $\{u+c\mathbf{1}|c\leq\frac{1}{\alpha}\log\frac{2\pi\chi(S)} {\sum_{i\in V} \overline{\mathcal{R}}_i e^{\alpha u_i}}\}$ stays in the set $\mathcal{A}$.
Hence, the set $\mathcal{A}$ is unbounded if $\alpha>0$.
This completes the proof.
\qed

According to Proposition \ref{Prop: ABC}, we have following result.
\begin{lemma}\label{Lem: minimum lies at the boundary}
If one of the following three conditions is satisfied
\begin{description}
  \item[(1)] $\alpha>0$ and the energy function $\mathcal{E}$ attains a minimum in the set $\mathcal{A}$,
  \item[(2)] $\alpha<0$ and the energy function $\mathcal{E}$ attains a minimum in the set $\mathcal{B}$,
  \item[(3)] $\alpha<0$ and the energy function $\mathcal{E}$ attains a minimum in the set $\mathcal{C}$,
\end{description}
then the minimum value point of $\mathcal{E}$ lies in the set $\{u\in \mathbb{R}^V|\sum_{i\in V} \overline{\mathcal{R}}_i e^{\alpha u_i}=2\pi\chi(S)\}$.
\end{lemma}
\proof
Suppose $\alpha>0$ and the function $\mathcal{E}$ attains a minimum at $u\in \mathcal{A}$.
The definition of $\mathcal{A}$ in (\ref{A}) implies $\chi(S)<0$.
Set
\begin{equation*}
c_0=\frac{1}{\alpha}\log\frac{2\pi\chi(S)} {\sum_{i\in V} \overline{\mathcal{R}}_i e^{\alpha u_i}},
\end{equation*}
then $c_0\geq0$.
By the proof of Proposition \ref{Prop: ABC}, $u+c_0\mathbf{1}\in \mathcal{A}$.
Therefore, by the additive property of the function $\mathcal{E}$ in (\ref{Eq: property of E}),
we have
\begin{equation*}
\mathcal{E}(u)\leq \mathcal{E}(u+c_0\mathbf{1})
=\mathcal{E}(u)+2\pi c_0\chi(S),
\end{equation*}
which implies $c_0\leq0$ by $\chi(S)<0$.
Hence $c_0=0$ and $\sum_{i\in V} \overline{\mathcal{R}}_i e^{\alpha u_i}=2\pi\chi(S)$.
This proves the case $\textbf{(1)}$.
The proofs for the cases $\textbf{(2)}$ and $\textbf{(3)}$ are similar, we omit the details here.
\qed

By Lemma \ref{Lem: minimum lies at the boundary}, we translate Theorem \ref{Thm: existence 3} into the following theorem, which is a non-convex optimization problem with inequality constraints.

\begin{theorem}\label{Thm: inequality constraints}
Let $(dist_S,r)$ be a decorated PE metric on a marked surface $(S,V)$ with $\chi(S)\neq0$.
Suppose $\alpha\in \mathbb{R}$ is a non-zero constant and $\overline{\mathcal{R}}$ is a given function defined on $V$.
\begin{description}
  \item[(1)] If $\overline{\mathcal{R}}\leq0$, $\overline{\mathcal{R}}\not\equiv0$, $\alpha>0$ and the energy function $\mathcal{E}$ attains a minimum in $\mathcal{A}$, then there exists a decorated PE metric in the discrete conformal class $\mathcal{D}(dist_S,r)$ with combinatorial $\alpha$-curvature $\overline{\mathcal{R}}$;
  \item[(2)]If $\overline{\mathcal{R}}>0$, $\alpha<0$ and the energy function $\mathcal{E}$ attains a minimum in $\mathcal{B}$, then there exists a decorated PE metric in the discrete  conformal class $\mathcal{D}(dist_S,r)$ with combinatorial $\alpha$-curvature $\overline{\mathcal{R}}$;
  \item[(3)] If $\overline{\mathcal{R}}\leq0$, $\overline{\mathcal{R}}\not\equiv0$, $\alpha<0$ and the energy function $\mathcal{E}$ attains a minimum in $\mathcal{C}$, then there exists a decorated PE metric in the discrete conformal class $\mathcal{D}(dist_S,r)$ with combinatorial $\alpha$-curvature $\overline{\mathcal{R}}$.
\end{description}
\end{theorem}
\proof
Lemma \ref{Lem: minimum lies at the boundary} shows that if $u\in \mathbb{R}^V$ is a minimum of the energy
function $\mathcal{E}$ defined on one of these sets, then $\sum_{i\in V} \overline{\mathcal{R}}_i e^{\alpha u_i}= 2\pi\chi(S)$.
The conclusion follows from the following claim.\\
\textbf{Claim :} Up to scaling, the decorated PE metrics with combinatorial $\alpha$-curvature $\overline{\mathcal{R}}$ in the discrete conformal class are in one-to-one correspondence with the critical points of the function $\mathcal{E}$ under the constraint $\sum_{i\in V} \overline{\mathcal{R}}_i e^{\alpha u_i}=2\pi\chi(S)$.

We use the method of Lagrange multipliers to prove this claim.
Set
\begin{equation*}
G(u,\mu)=\mathcal{E}(u)-\mu \left(\sum_{i\in V} \overline{\mathcal{R}}_i e^{\alpha u_i}-2\pi\chi(S)\right),
\end{equation*}
where $\mu\in \mathbb{R}$ is a Lagrange multiplier.
If $u$ is a critical point of the function $\mathcal{E}$ under the constraint $\sum_{i\in V} \overline{\mathcal{R}}_i e^{\alpha u_i}=2\pi\chi(S)$,
then by the fact $\nabla_{u_i} \mathcal{E}=K_i$, we have
\begin{equation*}
0=\frac{\partial G(u,\mu)}{\partial u_i}
=K_i-\mu\alpha\overline{\mathcal{R}}_i e^{\alpha u_i},
\end{equation*}
which implies
\begin{equation*}
\mathcal{R}_{\alpha, i}=\frac{K_i}{e^{\alpha u_i}}
=\mu\alpha\overline{\mathcal{R}}_i.
\end{equation*}
By the discrete Gauss-Bonnet formula (\ref{Eq: Gauss-Bonnet}), the Lagrange multiplier $\mu$ satisfies
\begin{equation*}
\mu=\frac{2\pi \chi(S)}{\alpha\sum_{i\in V} \overline{\mathcal{R}}_i e^{\alpha u_i}}=\frac{1}{\alpha}
\end{equation*}
under the constraint
$\sum_{i\in V} \overline{\mathcal{R}}_i e^{\alpha u_i}=2\pi\chi(S)$.
This implies the combinatorial $\alpha$-curvature
\begin{equation*}
\mathcal{R}_{\alpha, i}=\mu\alpha\overline{\mathcal{R}}_i
=\frac{2\pi \chi(S)}{\sum_{i\in V} \overline{\mathcal{R}}_i e^{\alpha u_i}}\overline{\mathcal{R}}_i=\overline{\mathcal{R}}_i
\end{equation*}
under the constraint
$\sum_{i\in V} \overline{\mathcal{R}}_i e^{\alpha u_i}=2\pi\chi(S)$.
\qed

\subsection{Reduction to Theorem \ref{Thm: main 2}}
By Theorem \ref{Thm: inequality constraints},
we just need to prove that the function $\mathcal{E}(u)$ attains the minimum in the sets $\mathcal{A},\ \mathcal{B}$ and $\mathcal{C}$ respectively.
Recall the following classical result from calculus.

\begin{theorem}\label{Thm: calculus}
Let $\Omega\subseteq \mathbb{R}^m$ be a closed set and $f: \Omega\rightarrow \mathbb{R}$ be a continuous function.
If every unbounded sequence $\{u_n\}_{n\in \mathbb{N}}$ in $\Omega$ has a subsequence $\{x_{n_k}\}_{k\in \mathbb{N}}$ such that $\lim_{k\rightarrow +\infty} f(x_{n_k})=+\infty$,
then $f$ attains a minimum in $\Omega$.
\end{theorem}

One can refer to \cite{Kourimska Thesis} (Section 4.1) for a proof of Theorem \ref{Thm: calculus}.
The majority of the conditions in Theorem \ref{Thm: calculus} are satisfied,
including the sets $\mathcal{A},\ \mathcal{B}$ and $\mathcal{C}$ are closed subsets of $\mathbb{R}^V$ by Proposition \ref{Prop: ABC} and the energy function $\mathcal{E}$ is continuous.
To prove Theorem \ref{Thm: existence 3},
we just need to prove the following theorem.

\begin{theorem}\label{Thm: main 2}
Suppose $(S,V)$ is a marked surface with a decorated PE metric $(dist_S,r)$, $\alpha\in \mathbb{R}$ is a constant and $\overline{\mathcal{R}}$ is a given function defined on $V$. If one of the following three conditions is satisfied
\begin{description}
  \item[(1)] $\alpha>0$ and $\{u_n\}_{n\in \mathbb{N}}$ is an unbounded sequence in $\mathcal{A}$,
  \item[(2)] $\alpha<0$ and $\{u_n\}_{n\in \mathbb{N}}$ is an unbounded sequence in $\mathcal{B}$,
  \item[(3)] $\alpha<0$ and $\{u_n\}_{n\in \mathbb{N}}$ is an unbounded sequence in $\mathcal{C}$,
\end{description}
then there exists a subsequence $\{u_{n_k}\}_{k\in \mathbb{N}}$ of $\{u_n\}_{n\in \mathbb{N}}$ such that
$\lim_{k\rightarrow +\infty} \mathcal{E}(u_{n_k})=+\infty$.
\end{theorem}

\subsection{Behaviour of sequences of conformal factors}
Let $\{u_n\}_{n\in \mathbb{N}}$ be an unbounded sequence in $\mathbb{R}^V$.
Denote its coordinate sequence at $j\in V$ by $\{u_{j,n}\}_{n\in \mathbb{N}}$.
Motivated by \cite{Kourimska}, we call the sequence $\{u_n\}_{n\in \mathbb{N}}$ with the following properties as a ``good" sequence.
\begin{description}
\item[(1)] It lies in one cell $\mathcal{C}_\mathcal{T}(dist_{S},r)$ of $\mathbb{R}^V$;
\item[(2)] There exists a vertex $i^*\in V$ such that $u_{i^*,n}\leq u_{j,n}$ for all $j\in V$ and $n\in \mathbb{N}$;
\item[(3)] Each coordinate sequence $\{u_{j,n}\}_{n\in \mathbb{N}}$ either converges, diverges properly to $+\infty$, or diverges properly to $-\infty$;
\item[(4)] For any $j\in V$, the sequence $\{u_{j,n}-u_{i^*,n}\}_{n\in \mathbb{N}}$ either converges or diverges properly to $+\infty$.
\end{description}

By Lemma \ref{Lem: finite decomposition},
it is obvious that every sequence of discrete conformal factors in $\mathbb{R}^V$ possesses a ``good" subsequence.
Hence, the ``good" sequence could be chosen without loss of generality.
To prove Theorem \ref{Thm: main 2}, we further need the following two results obtained by the authors in \cite{XZ}.

\begin{lemma}[\cite{XZ}, Corollary 3.6]\label{Lem: one infty two converge}
For a discrete conformal factor $u\in \mathbb{R}^V$, let $\mathcal{T}$ be a weighted Delaunay triangulation of the decorated PE surface $(S,V,dist_S(u),r(u))$.
For any decorated triangle $\{ijk\}\in F$ in $\mathcal{T}$,
at least two of the three sequences $\{u_{i,n}-u_{i^*,n}\}_{n\in \mathbb{N}}$,
$\{u_{j,n}-u_{i^*,n}\}_{n\in \mathbb{N}}$,
$\{u_{k,n}-u_{i^*,n}\}_{n\in \mathbb{N}}$
converge.
\end{lemma}

\begin{lemma}[\cite{XZ}, Lemma 3.12]\label{Lem: E decomposition}
There exists a convergent sequence $\{D_n\}_{n\in \mathbb{N}}$ such that the function $\mathcal{E}$ satisfies
\begin{equation*}
\mathcal{E}(u_n)=D_n+2\pi\left(u_{i^*,n}\chi(S)+\sum_{j\in V}(u_{j,n}-u_{i^*,n})\right).
\end{equation*}
\end{lemma}

\noindent\textbf{Proof of Theorem \ref{Thm: main 2}:}
Let $\{u_n\}_{n\in \mathbb{N}}$ be an unbounded ``good" sequence.
We just need to prove that $\lim_{n\rightarrow +\infty} \mathcal{E}(u_n)=+\infty$.

\noindent\textbf{(1):}
Let $\alpha>0$ and $\{u_n\}_{n\in \mathbb{N}}$ be an unbounded sequence in $\mathcal{A}$.
The definition of $\mathcal{A}$ in (\ref{A}) implies $\chi(S)<0$, $\overline{\mathcal{R}}\leq0$ and $\overline{\mathcal{R}}\not\equiv0$.
Since the sequence $\{u_n\}_{n\in \mathbb{N}}$ lies in $\mathcal{A}$, we have
\begin{equation}\label{Eq: key 1}
0>\sum_{j\in V} \overline{\mathcal{R}}_j e^{\alpha (u_{j,n}-u_{i^*,n})}
=e^{-\alpha u_{i^*,n}}\cdot\sum_{j\in V} \overline{\mathcal{R}}_j e^{\alpha u_{j,n}}
\geq 2\pi\chi(S) e^{-\alpha u_{i^*,n}}.
\end{equation}
By the definition of ``good" sequence, the sequence $\left\{\sum_{j\in V}(u_{j,n}-u_{i^*,n})\right\}_{n\in \mathbb{N}}$ converges to a finite positive number or diverges properly to $+\infty$

If $\left\{\sum_{j\in V}(u_{j,n}-u_{i^*,n})\right\}_{n\in \mathbb{N}}$ converges to a finite positive number,
then the sequence $\{u_{j,n}-u_{i^*,n}\}_{n\in \mathbb{N}}$ converges for any $j\in V$.
This implies $\sum_{j\in V} \overline{\mathcal{R}}_j e^{\alpha (u_{j,n}-u_{i^*,n})}$ converges to a finite negative number by (\ref{Eq: key 1}).
Then by $\chi(S)<0$, we have
\begin{equation*}
-\alpha u_{i^*,n}\geq \ln\frac{\sum_{j\in V} \overline{\mathcal{R}}_j e^{\alpha (u_{j,n}-u_{i^*,n})}}{2\pi\chi(S)}.
\end{equation*}
Hence $\{u_{i^*,n}\}_{n\in \mathbb{N}}$ is bounded from above by $\alpha>0$.
This implies $\{u_{i^*,n}\}_{n\in \mathbb{N}}$ converges to a finite number or diverges properly to $-\infty$.
If $\{u_{i^*,n}\}_{n\in \mathbb{N}}$ converges to a finite number,
then by $\{u_{j,n}-u_{i^*,n}\}_{n\in \mathbb{N}}$ converges for any $j\in V$,
we have $\{u_{j,n}\}_{n\in \mathbb{N}}$ is bounded for any $j\in V$.
This contradicts the assumption that $\{u_n\}_{n\in \mathbb{N}}$ is unbounded.
Therefore, the sequence $\{u_{i^*,n}\}_{n\in \mathbb{N}}$ diverges properly to $-\infty$.
Combining this with $\chi(S)<0$ and Lemma \ref{Lem: E decomposition},
we have $\lim_{n\rightarrow +\infty} \mathcal{E}(u_n)=+\infty$.

If $\left\{\sum_{j\in V}(u_{j,n}-u_{i^*,n})\right\}_{n\in \mathbb{N}}$ diverges properly to $+\infty$,
then there exists at least one vertex $j\in V$ such that the sequence $\{u_{j,n}-u_{i^*,n}\}_{n\in \mathbb{N}}$ diverges properly to $+\infty$.
By Lemma \ref{Lem: one infty two converge},
for any vertex $k\sim j$, the sequence $\{u_{k,n}-u_{i^*,n}\}_{n\in \mathbb{N}}$ converges.
Since $\alpha>0$, then $e^{\alpha (u_{j,n}-u_{i^*,n})}$ converges to a finite positive number or diverges properly to $+\infty$ and
for at least one vertex $j\in V$ the term $e^{\alpha (u_{j,n}-u_{i^*,n})}$ converges to a finite positive number.
Since $\overline{\mathcal{R}}\leq0$ and $\overline{\mathcal{R}}\not\equiv0$,
then
$\sum_{j\in V} \overline{\mathcal{R}}_j e^{\alpha (u_{j,n}-u_{i^*,n})}$
converges to a finite negative number or diverges properly to $-\infty$.
\begin{description}
\item[$(i)$]
Suppose $\sum_{j\in V} \overline{\mathcal{R}}_j e^{\alpha (u_{j,n}-u_{i^*,n})}$ converges to a finite negative number.
Similar arguments imply $u_{i^*,n}$ is bounded from above, then $u_{i^*,n}\chi(S)$ is bounded from below by $\chi(S)<0$.
Combining with the assumption that
$\left\{\sum_{j\in V}(u_{j,n}-u_{i^*,n})\right\}_{n\in \mathbb{N}}$ diverges properly to $+\infty$, we have
$\lim_{n\rightarrow +\infty} \mathcal{E}(u_n)=+\infty$ by Lemma \ref{Lem: E decomposition}.
\item[$(ii)$]
Suppose $\sum_{j\in V}\overline{\mathcal{R}}_j e^{\alpha (u_{j,n}-u_{i^*,n})}$ diverges properly to $-\infty$.
Then $2\pi\chi(S) e^{-\alpha u_{i^*,n}}$ diverges properly to $-\infty$ by (\ref{Eq: key 1}).
Since $\chi(S)<0$, then $e^{-\alpha u_{i^*,n}}$ diverges properly to $+\infty$.
By $\alpha>0$, then $\{u_{i^*,n}\}_{n\in \mathbb{N}}$ diverges properly to $-\infty$.
Hence $u_{i^*,n}\chi(S)$ diverges properly to $+\infty$ by $\chi(S)<0$.
Then $\lim_{n\rightarrow +\infty} \mathcal{E}(u_n)=+\infty$ by Lemma \ref{Lem: E decomposition}.
\end{description}

\noindent\textbf{(2):} Let $\alpha<0$ and $\{u_n\}_{n\in \mathbb{N}}$ be an unbounded sequence in $\mathcal{B}$. The definition of $\mathcal{B}$ in (\ref{B}) implies $\chi(S)>0$ and $\overline{\mathcal{R}}>0$.
Since the sequence $\{u_n\}_{n\in \mathbb{N}}$ lies in $\mathcal{B}$, we have
\begin{equation}\label{Eq: key 2}
0<\sum_{j\in V} \overline{\mathcal{R}}_j e^{\alpha (u_{j,n}-u_{i^*,n})}
=e^{-\alpha u_{i^*,n}}\cdot\sum_{j\in V} \overline{\mathcal{R}}_j e^{\alpha u_{j,n}}
\leq 2\pi\chi(S) e^{-\alpha u_{i^*,n}}.
\end{equation}

If $\left\{\sum_{j\in V}(u_{j,n}-u_{i^*,n})\right\}_{n\in \mathbb{N}}$ converges,
then the sequence $\{u_{j,n}-u_{i^*,n}\}_{n\in \mathbb{N}}$ converges for any $j\in V$.
This implies $\sum_{j\in V} \overline{\mathcal{R}}_j e^{\alpha (u_{j,n}-u_{i^*,n})}$ converges to a finite positive number by (\ref{Eq: key 2}).
Since $\chi(S)>0$, then the equation (\ref{Eq: key 2}) implies
\begin{equation*}
-\alpha u_{i^*,n}\geq \ln \frac{\sum_{j\in V} \overline{\mathcal{R}}_j e^{\alpha (u_{j,n}-u_{i^*,n})}}{2\pi\chi(S)}.
\end{equation*}
By $\alpha<0$, then $\{u_{i^*,n}\}_{n\in \mathbb{N}}$ is bounded from below.
This implies $\{u_{i^*,n}\}_{n\in \mathbb{N}}$ converges to a finite number or diverges properly to $+\infty$.
Combining this with $\{u_n\}_{n\in \mathbb{N}}$ is unbounded
and $\{u_{j,n}-u_{i^*,n}\}_{n\in \mathbb{N}}$ converges for all $j\in V$,
we have the sequence $\{u_{i^*,n}\}_{n\in \mathbb{N}}$ diverges properly to $+\infty$.
By $\chi(S)>0$ and Lemma \ref{Lem: E decomposition},
we have $\lim_{n\rightarrow +\infty} \mathcal{E}(u_n)=+\infty$.

If the sequence $\left\{\sum_{j\in V}(u_{j,n}-u_{i^*,n})\right\}_{n\in \mathbb{N}}$ diverges properly to $+\infty$,
then there exists at least one vertex $j\in V$ such that the sequence $\{u_{j,n}-u_{i^*,n}\}_{n\in \mathbb{N}}$ diverges properly to $+\infty$.
By Lemma \ref{Lem: one infty two converge},
for any vertex $k\sim j$, the sequence $(u_{k,n}-u_{i^*,n})_{n\in \mathbb{N}}$ converges.
Therefore, $e^{\alpha (u_{j,n}-u_{i^*,n})}$ converges to zero or a finite positive number and
for at least one vertex $j\in V$ the term $e^{\alpha (u_{j,n}-u_{i^*,n})}$ converges to a finite positive number.
Since $\alpha<0$ and $\overline{\mathcal{R}}>0$,
then $\sum_{j\in V} \overline{\mathcal{R}}_j e^{\alpha (u_{j,n}-u_{i^*,n})}$ converges to a finite positive number.
This implies
$2\pi\chi(S) e^{-\alpha u_{i^*,n}}$ has a positive lower bound by (\ref{Eq: key 2}).
By $\alpha<0$ and $\chi(S)>0$, then $\{u_{i^*,n}\}_{n\in \mathbb{N}}$ is bounded from below.
Then $u_{i^*,n}\chi(S)$ is bounded from below.
Combining this with $\left\{\sum_{j\in V}(u_{j,n}-u_{i^*,n})\right\}_{n\in \mathbb{N}}$ diverges properly to $+\infty$,
we have $\lim_{n\rightarrow +\infty} \mathcal{E}(u_n)=+\infty$ by Lemma \ref{Lem: E decomposition}.

\noindent\textbf{(3):}
Let $\alpha<0$ and $\{u_n\}_{n\in \mathbb{N}}$ be an unbounded sequence in $\mathcal{C}$.
The definition of $\mathcal{C}$ in ($\ref{C}$) implies $\chi(S)<0$ and $\overline{\mathcal{R}}\leq0$ and $\overline{\mathcal{R}}\not\equiv0$.
Since the sequence $\{u_n\}_{n\in \mathbb{N}}$ lies in $\mathcal{C}$, we have
\begin{equation}\label{Eq: key 3}
\sum_{j\in V} \overline{\mathcal{R}}_j e^{\alpha (u_{j,n}-u_{i^*,n})}=e^{-\alpha u_{i^*,n}}\cdot\sum_{j\in V} \overline{\mathcal{R}}_j e^{\alpha u_{j,n}}\leq 2\pi\chi(S)e^{-\alpha u_{i^*,n}}<0.
\end{equation}

If $\left\{\sum_{j\in V}(u_{j,n}-u_{i^*,n})\right\}_{n\in \mathbb{N}}$ converges, then the sequence $\{u_{j,n}-u_{i^*,n}\}_{n\in \mathbb{N}}$ converges for all $j\in V$.
This implies that
$\sum_{j\in V} \overline{\mathcal{R}}_j e^{\alpha (u_{j,n}-u_{i^*,n})}$ converges to a finite negative number by (\ref{Eq: key 3}).
Since $\chi(S)<0$, then the equation (\ref{Eq: key 3}) implies
\begin{equation*}
-\alpha u_{i^*,n}\leq \ln \frac{\sum_{j\in V} \overline{\mathcal{R}}_j e^{\alpha (u_{j,n}-u_{i^*,n})}}{2\pi\chi(S)}.
\end{equation*}
By $\alpha<0$, we have $\{u_{i^*,n}\}_{n\in \mathbb{N}}$ is bounded from above.
Combining this with $\{u_n\}_{n\in \mathbb{N}}$ is unbounded
and $\{u_{j,n}-u_{i^*,n}\}_{n\in \mathbb{N}}$ converges for all $j\in V$,
the sequence $\{u_{i^*,n}\}_{n\in \mathbb{N}}$ diverges properly to $-\infty$.
Combining this with $\chi(S)<0$ and Lemma \ref{Lem: E decomposition},
we have $\lim_{n\rightarrow +\infty} \mathcal{E}(u_n)=+\infty$.

If $\left\{\sum_{j\in V}(u_{j,n}-u_{i^*,n})\right\}_{n\in \mathbb{N}}$ diverges properly to $+\infty$,
then there exists at least one vertex $j\in V$
such that the sequence $\{u_{j,n}-u_{i^*,n}\}_{n\in \mathbb{N}}$ diverges properly to $+\infty$.
By Lemma \ref{Lem: one infty two converge},
for any vertex $k\sim j$, the sequence $\{u_{k,n}-u_{i^*,n}\}_{n\in \mathbb{N}}$ converges.
Since $\alpha<0$, then $e^{\alpha (u_{j,n}-u_{i^*,n})}$ converges to zero or a finite positive number and
for at least one vertex $j\in V$ the term $e^{\alpha (u_{j,n}-u_{i^*,n})}$ converges to a finite positive number.
Note that $\sum_{j\in V} \overline{\mathcal{R}}_j e^{\alpha (u_{j,n}-u_{i^*,n})}<0$ by (\ref{Eq: key 3}).
Therefore, $\sum_{j\in V} \overline{\mathcal{R}}_j e^{\alpha (u_{j,n}-u_{i^*,n})}$ converges to zero or a finite negative number.
\begin{description}
\item[$(i)$]
Suppose $\sum_{j\in V} \overline{\mathcal{R}}_j e^{\alpha (u_{j,n}-u_{i^*,n})}$ converges to zero.
Then $2\pi\chi(S) e^{-\alpha u_{i^*,n}}$ converges to zero by (\ref{Eq: key 3}).
Since $\alpha<0$ and $\chi(S)<0$, then $\{u_{i^*,n}\}_{n\in \mathbb{N}}$ diverges properly to $-\infty$.
Combining this with $\chi(S)<0$ and Lemma \ref{Lem: E decomposition},
we have $\lim_{n\rightarrow +\infty} \mathcal{E}(u_n)=+\infty$.

\item[$(ii)$]
Suppose $\sum_{j\in V} \overline{\mathcal{R}}_j e^{\alpha (u_{j,n}-u_{i^*,n})}$ converges to a finite negative number.
Then $2\pi\chi(S) e^{-\alpha u_{i^*,n}}$ has a negative lower bound by (\ref{Eq: key 3}).
By $\alpha<0$ and $\chi(S)<0$, then $u_{i^*,n}$ is bounded from above.
Combining this with $\chi(S)<0$ and $\left\{\sum_{j\in V}(u_{j,n}-u_{i^*,n})\right\}_{n\in \mathbb{N}}$
diverges properly to $+\infty$,
we have $\lim_{n\rightarrow +\infty} \mathcal{E}(u_n)=+\infty$ by Lemma \ref{Lem: E decomposition}.
\end{description}
\qed

\end{document}